\documentclass[11pt]{article}
\usepackage[dvips]{graphics,color}
\def\hind{\hangindent=2pc\hangafter=1}
\newfont{\smcaps}{cmcsc10 scaled\magstep1}
\newfont{\subsubti}{cmss12 scaled\magstep1}
\newfont{\tti}{cmssi10 scaled\magstep1}
\newfont{\subti}{cmssbx10 scaled\magstep2}
\newfont{\ti}{cmssbx10 scaled\magstep3}
\newfont{\nb}{cmssi10 scaled\magstep1}
 at 11truept

\newcommand{\IMA}{${{\rm ~IMA}\,}$}
\newcommand{\AR}{{\rm ~AR\,}}
\newcommand{\ARIMA}{{\rm ~ARIMA\,}}

\newcommand{\SE}{{\rm ~SE\,}}
\def\Var{\rm ~Var\,}
\def\Pr{\rm ~Pr\,}

\raggedright
\begin{document}

\title{\ti Power Computations for Intervention Analysis}
\date{}
\author{A. Ian McLeod and Evelyn R. Vingilis\\
Department of Statistical and Actuarial Sciences,\\
The University of Western Ontario\\ London, Ontario N6A 5B7 \\
aimcleod@uwo.ca\\\\
Department of Family Medicine,\\ University of Western Ontario \\
London, Ontario N6G 4X8 \\
evingili@uwo.ca \\\\
}

\maketitle

\hrule

\bigskip

\noindent A. Ian McLeod and Evelyn R. Vingilis (2005),
Power Computations for Intervention Analysis,
{\it Technometrics} 47/2, 174-180.

\newpage

\centerline{Abstract}
\medskip

In many intervention analysis applications time series data may be expensive
or otherwise difficult to collect.
In this case the power function is helpful since it
can be used to determine the probability that a proposed intervention
analysis application will detect a meaningful change.
Assuming that an underlying ARIMA or fractional ARIMA model
is known or can be estimated from the pre-intervention time series,
the methodology for computing the required
power function is developed for pulse, step and ramp interventions
with ARIMA and fractional ARIMA errors.
Convenient formulae for computing the power function
for important special cases are given.
Illustrative applications in traffic safety and environmental impact assessment
are discussed.


\newpage

\centerline{\subti 1. INTRODUCTION}
\medskip

Intervention analysis developed by Box and Tiao (1976a) has been
widely used in a variety of applications in engineering, biological,
environmental and social sciences to quantify the
effect of a known intervention at time $t=T$ on data collected
as a time series, $z_t,\ t=1,\ldots,n$.
In its simplest form,
intervention analysis itself may be regarded as a generalization
of the two-sample problem to the case where the error or noise term is
autocorrelated.
It is well-known that the usual two-sample procedures are not robust against alternatives
involving autocorrelation (Box, Hunter and Hunter, 1978, \S 3.1).
The purpose of this article is to describe methods for computing the necessary
sample size to detect an intervention with a prescribed power and level.
It is shown by simulation experiments that these methods can be accurate even
in moderately small samples.
Statistical power computations have also been studied by
Tiao et al. (1990) and Weatherhead et al. (1998) for particular types of intervention
analysis models used for trend detection with environmental time series.
This article extends and refines these results.

It is assumed that for $t<T+b$, where $b$ is the delay parameter,
the time series is generated by
a fractional $\ARIMA(p,d,q)$ with fractional differencing parameter $|f|<0.5$.
Stationary short-memory time series models, $d=f=0$, are used in
environmental impact assessment
(Box and Tiao, 1976a; Tiao et al., 1990; Noakes and Campbell, 1992; Weatherhead et al. 1998;  Hipel and McLeod, 1994, \S 19.4.5)
and in quality control (Jiang, Tsui and Woodall, 2000)
as well as in many other areas of science and technology.
Nonstationary models with $d=1$ and/or long-memory models with $0<f<0.5$
have numerous applications in the physical and engineering sciences
such as:
quality control and industrial time series (Luce\~no, 1995; Box and Luce\~no, 1997),
internet traffic (Cao et al., 2001),
daily solar irradiance (K\"arner, 2002),
levels of Lake Huron (Roberts, 1991, p.319-320),
daily wind-speed (Haslett and Raftery, 1989),
and various types of hydrological time series
(Beran, 1994; Hipel and McLeod, 1994).

In general, we may write the fractional $\ARIMA$ model for the pre-intervention series as
\begin{equation}
\nabla^{d+f}  z_t =
\xi + \theta(B) / \phi(B) a_t,\quad t=1,\ldots,T+b-1,
\label{ARIMAModel}
\end{equation}
where $\xi$ is the constant term,
$d$ is the differencing parameter,
$\nabla=1-B$,
$\theta(B)=1-\theta_1 B - \ldots -\theta_q B^q$,
$\phi(B)=1-\phi_1 B - \ldots -\phi_p B^p$
and
$B$ is the backshift operator on $t$.
The innovations, denoted by $a_t,\ t=1,\ldots,n$,
are assumed to be independent and normally distributed with mean zero and variance $\sigma_a^2$.
It is also assumed that $\phi(B)=0$ and $\theta(B)=0$ have no common roots and that all
roots are outside the unit circle.
\bigskip

\centerline{\subti 2. SIMPLE INTERVENTION ANALYSIS (SIA) MODEL}

\medskip
{\subsubti\noindent 2.1 Introduction}
\smallskip

The SIA model may be written,
\begin{equation}
\nabla^d  z_t =
\xi+
\omega \nabla^d B^b I_t^{(T)} +
\nabla^{-f} \frac{\theta(B)}{\phi(B)} a_t,\quad t=1,\ldots,n,
\label{SimpleStepInterventionModel}
\end{equation}
where $I_t^{(T)}$ is the intervention series,
$\omega$ is the parameter indicating the magnitude of the intervention
and
$\nabla^{-f} \theta(B) / \phi(B) a_t$ is the stationary error component.
In this article three types of intervention series are used,
the step, pulse and ramp series, defined respectively by,
\begin{equation}
I_t^{(T)}= S_t^{(T)}= \cases{0, &if $t<T$,\cr 1& if $t\ge T$,\cr}
\label{PulseIA}
\end{equation}
\begin{equation}
I_t^{(T)}=P_t^{(T)}= \cases{0, &if $t \ne T$,\cr 1& if $t = T$,\cr}
\label{StepIA}
\end{equation}
or
\begin{equation}
I_t^{(T)}=R_t^{(T)}= \cases{0, &if $t<T$,\cr t-T+1 & if $t\ge T$.\cr}
\label{RampIA}
\end{equation}
In practice two of the most common models for the error are
the $\AR(1)$ and $\IMA(1)$ which correspond
respectively to $p=1, d=0, q=0$ and $p=0, d=1, q=1$.
In the case of a step intervention, the SIA model implies that for
$t\ge T+b$ an increase of $\omega$ occurred.
So the SIA model with a step intervention can be regarded as the time-series generalization of the standard two-sample test for a change in location and in practice this is one
of the most frequently applicable models.
Pulse interventions are useful for dealing with outliers (Chang, Tiao and Chen, 1988).
A ramp intervention has been used to model the recovery trend in stratospheric ozone
(Reinsel et al. 2002).

The SIA model may be generalized by allowing for multiple interventions and other types of interventions, as well as for seasonal ARIMA errors and possible covariates
(Tiao et al., 1990; Weatherhead et al., 1998; Reinsel, 2002; Reinsel et al., 2002).
All of these situations are easily handled
with the methods discussed in \S 1.2 and \S 1.3.
Power computations, although possible, are less useful when applied to dynamic response
interventions for the reasons explained in Appendix B.

\medskip
{\subsubti\noindent 2.2 Information Matrix}
\smallskip

Letting $\lambda_1=(\xi, \omega)$ and $\lambda_2=(\phi_1,\ldots,\phi_p,\theta_1,\ldots,\theta_q,f)$,
it is shown in Appendix A that the expected Fisher information matrix
is block diagonal with blocks, ${\cal I}_{\lambda_1}$ and ${\cal I}_{\lambda_2}$ corresponding to $\lambda_1$ and $\lambda_2$.
For the first block,
\begin{equation}
{\cal I}_{\xi,\omega}=\sigma_a^{-2} J^{\prime} \Gamma_n^{-1} J,
\label{InformationMatrixGeneral}
\end{equation}
where
$\sigma_a^{-2}  \Gamma_n^{-1}$ is the inverse of the covariance matrix
of the stationary component and
$J$ is an $n \times 2$ matrix with $1$ in the first column and
$\nabla^d I_t^{(T)}, t=1,\dots,n$ in the second column.
The Trench algorithm (Golub and Van Loan, 1983) provides a computationally
efficient method for computing $\Gamma_n^{-1}$.
An expression essentially equivalent to eqn. (\ref{InformationMatrixGeneral})
was obtained by Tiao et al. (1990) and Weatherhead et al. (1998) using generalized least squares.
Assuming approximate normality of the estimates, the asymptotic variance of
the maximum likelihood estimate of $\omega$ is found by taking
the $(2,2)$ element of the inverse of (\ref{InformationMatrixGeneral}),
\begin{equation}
\sigma_{\hat \omega} = \surd \left( {\cal I}_{1,1}/\left({\cal I}_{1,1} {\cal I}_{2,2} - {\cal I}_{1,2}^2\right) \right),
\label{SEOmega}
\end{equation}
where ${\cal I}_{i,j}$ denotes the $(i,j)$ entry in the matrix ${\cal I}_{\xi,\omega}$.
If the constant term, $\xi$, is not present,
$\sigma_{\hat \omega} = 1/\surd {\cal I}_{2,2}$.
When there is an extensive amount of data prior to the intervention
it is sometimes helpful to simply correct the series
by its sample mean and assume $\xi = 0$ (Tiao et al., 1990).

The results of Pierce (1972) provide a computationally efficient approximation to (\ref{InformationMatrixGeneral}) when $f=0$.
From Pierce (1972, eqn. 3.2)
we can write the Fisher information for $(\xi,\omega)$ based on $n$ observations as
\begin{equation}
{\cal I}_{\xi,\omega}=\sigma_a^{-2} \pmatrix{n\kappa^2& \kappa \sum_t v_t&\cr
\kappa \sum_t v_t&  \sum_t v_t^2&\cr},
\label{InformationMatrixPierce}
\end{equation}
where
$\kappa=-\phi(1)/\theta(1)$ and $v_t = -\phi(B)/\theta(B)w_t$,
where $w_t = \nabla^d I_t^{(T)}$.  Without loss of generality
we take $b=0$ since if $b>0$, the formulae hold with $T$ replaced
by $T+b$.
Provided that $T$ is not too small and $T$ is not too close to $n$, eqn. (\ref{InformationMatrixPierce}) yields almost identical values to the more
exact formula given in (\ref{InformationMatrixGeneral}).
New explicit expressions, using Pierce's approximation for $\AR(1)$ and $\IMA(1)$ cases, are given in
Tables 1 and 2 below for step, pulse and ramp interventions.

\centerline{[Tables 1 and 2 about here]}

From eqn. (\ref{InformationMatrixGeneral}), it follows that for consistency
of the estimates $\hat \xi$ and $\hat\omega$, ${\cal I}_{\xi,\omega}/n$
or equivalently, $J^{\prime} J/n$, must converge to a nonsingular matrix.
For the intervention analysis models defined by eqns.
(\ref{SimpleStepInterventionModel}), (\ref{PulseIA}), (\ref{StepIA}) and (\ref{RampIA}),
this happens provided that
\begin{equation}
{1\over n}\sum\displaylimits_{t=1}^n \nabla^d I_t^{(T)} \rightarrow c,\quad c>0, c\ne 1.
\label{PierceCondition}
\end{equation}
If the constant term, $\xi$, is assumed to be known or zero then only $c>0$ is needed.
This result is certainly not the whole story from the application point of view.
In \S 1.5 we show using simulation experiments that the empirical variances
may be accurately estimates from \ref{SEOmega}) even when
eqn. ({\ref{PierceCondition}) is not satisfied.

\medskip
{\subsubti\noindent 2.3 Power and Sample Size}
\smallskip

The null hypothesis ${\cal H}_0:\omega=0$ can be tested using two asymptotically
equivalent methods.
The first method, referred to as the $Z$-test, uses
$Z=\hat \omega/{\hat\sigma_{\hat \omega}}$,
where  $\hat \omega$ is the maximum likelihood estimate for $\omega$
and ${\hat\sigma_{\hat \omega}}$ is its estimated standard error.
Note that ${\sigma_{\hat \omega}}$, the standard error of $\hat \omega$, depends only
on the underlying ARIMA model in the pre-intervention period and so it
can be estimated before the post-intervention data are obtained.
A second asymptotically equivalent method is to use a likelihood-ratio test.

The asymptotic theoretical power function for the $Z$-test of the null hypothesis
${\cal H}_0:\omega=0$
against the two-sided alternative at level $\alpha$ is
$\Pr\{|\ \hat \omega \ |>{\cal Z}_{1-\alpha/2}\ {\sigma_{\hat \omega}} | \omega \}$,
where ${\cal Z}_{1-\alpha/2}$ is the upper $(1-\alpha/2)$-quantile in the standard normal distribution.
For brevity the asymptotic theoretical power function will be referred to
simply as the power function.
In practice this power function is approximated by replacing $\sigma_{\hat \omega}$
by an estimate, $\hat \sigma_{\hat \omega}$, based either on the pre-intervention data
or on other prior knowledge.
Often it is more convenient to use the rescaled parameter,
$\delta=\omega/\sigma$, where
$\sigma^2$ is the variance of the stationary error component since in this
case knowledge of $\sigma^2$ is not needed.
The power function may be expressed in terms of $\delta$ as
\begin{equation}
\Pi(\delta)=
\Phi(-{\cal Z}_{1-\alpha/2}-\delta \sigma / {\sigma_{\hat\omega}})+
1-\Phi({\cal Z}_{1-\alpha/2}-\delta \sigma / {\sigma_{\hat\omega}}),
\label{PowerFunctionTwoSided}
\end{equation}
where $\Phi(\bullet)$ denotes the cumulative distribution function of the standard normal.
If the variance of the pre-intervention series, $\sigma^2$, is known or estimated,
the power function for $\omega$ is $\Pi(\omega/\sigma)$.
Eqn. (\ref{PowerFunctionTwoSided}) should be adjusted if only
a one-sided alternative is under consideration.

As in Tiao et al. (1990) it is sometimes of interest to estimate
the amount of additional data needed to detect an intervention of a specified
magnitude with a prescribed power.
The power function $\Pi(\delta)$ may be expressed more fully as a function of
the test level $\alpha$ and
the other underlying parameters $n$ and $T$ so we can write
the power function more fully as $\Pi(\delta,\alpha, n,T)$.
For a fixed $\alpha=\alpha^{(0)}$, $\delta=\delta^{(0)}$
and a prescribed power $\Pi^{(0)}$ we may estimate the number of
additional data values, $m$, that are required by numerically solving the equation
$\Pi(\delta^{(0)},\alpha^{(0)}, T+m-1 ,T)$\ $=\Pi^{(0)}$.
If as in the geophysical datasets considered in Tiao et al. (1990) there is
extensive pre-intervention data, we may assume the mean is known and take $T=1$ and solve $\Pi(\delta^{(0)},\alpha^{(0)}, m ,1)$\ $=\Pi^{(0)}$.
This technique is illustrated in \S 1.4 where it is also explained that in some
situations, due to the limitations imposed by the model, there is no solution for $m$.

In general the power and sample size computations for interventions
with ARIMA and fractional ARIMA errors
are easily done using an advanced {\it quantitative programming environment\/} such as
{\it Mathematica\/}, MatLab, S or Stata.
In the case of SIA with $\AR(1)$ or $\IMA(1)$ errors, power computations can
even be done on a hand calculator.

\medskip
{\subsubti\noindent 2.4 Numerical Illustrations}
\smallskip

\newcounter{AsymptoticPower}
\setcounter{AsymptoticPower}{1}
\newcounter{FigureSimulatedPower}
\setcounter{FigureSimulatedPower}{2}
\newcounter{FigureQPower}
\setcounter{FigureQPower}{3}

The power and sample size computations are illustrated in this section for the
SIA with a step intervention with $\AR(1)$, $\IMA(1)$ and fractionally-differenced
white noise.
First an approximation to the detection limit, $\delta^\prime$, is derived
for the step intervention in an SIA model with unknown mean, stationary short-memory
errors, with $f=d=0$, and a fixed number, $T-1$, of pre-intervention observations.
The variance of the estimate, $\hat \delta$, may be written,
$\Var(\hat \delta) \doteq \gamma_{\delta}/T$, where
$\gamma_{\delta} =$\ $  \sum_{k=-\infty}^{\infty} \gamma_k / \gamma_0$,
$\gamma_k$ is the autocovariance function for the stationary pre-intervention series
and $\gamma_0 = \sigma^2$.
To achieve 90\% power, $\Pr\{(\hat \delta-\delta^\prime)/\SE(\hat \delta)>1.96$
$-\delta^\prime/\SE(\hat \delta)\}$ $ \doteq 0.9.$
Hence $2-\delta^\prime/\SE(\hat \delta) \doteq -1.3$.
So $\delta^\prime \doteq 3.3\SE(\hat \delta)$.

Using Table 1, the power curve
for the $\AR(1)$ with unknown mean, $n=50$, $T=25$ and $\phi_1=0.5$,
$\sigma_{\omega}=0.526681$.
With $\sigma= $\ $1/\surd (1-\phi_1^2) =$\ $ 1.1547$, the power curve is
$\Pi(\delta)=$
$1+\Phi(-1.960 - 2.192\times \delta)$
$-\ \Phi(1.960 - 2.192\times \delta)$.
This and the power curve obtained by letting $n \rightarrow \infty$ are
shown in Figure \theAsymptoticPower\
as well as the approximate detection level, $\delta^{\prime}\doteq \gamma_{\delta}/\surd T=1.14$.
For comparison, the exact value of $\delta^{\prime}$ found by numerically solving
$\Pi(\delta^{\prime},0.5, 10^9 ,25)=0.9$
is
$\delta^{\prime}=1.12$.
Assuming an unknown mean and that $T=25$, we can find $m$, the number of additional
observations needed to achieve a prescribed power level.
For example, for 90\% power with $\delta^{(0)}=1.5$, solving
$\Pi(1.5,0.05, 25+m-1 ,25)=0.9$ we find $m=23$.
In the known mean case taking $T=1$ we find $m=10$.
In the unknown mean case, if $\delta^{(0)} \le \gamma_{\delta}$ there is no solution
but if the mean is known then $m$ can always be found.

\medskip
\centerline{[Figure \theAsymptoticPower\ about here]}
\medskip

The middle panels of Figure  \theFigureSimulatedPower\
illustrate the power curves for an $\IMA(1)$ with $n=50$ and $T=25$.
With $\theta_1=0.5$,
$\Pi(\delta)=$
$1+\Phi(-1.960 - 1.252\times \delta)$
$-\Phi(1.960 - 1.252\times \delta)$.

Since long-memory or fractional time series have also been suggested for
various types of geophysical data, it is of interest to examine the impact
of this type of process on our ability to detect interventions.
Table 3 compares the power of a two-sided 5\% level test of
the fractionally differenced white noise model $p=d=q=0$ with $f=0.2$
and $f=0.4$ to the corresponding approximating ARMA$(1,1)$ when $n=50$ and $T=25$.
The approximating ARMA$(1,1)$ model was determined by equating the first two
autocorrelations in the fractional model with the first two autocorrelations
in the ARMA$(1,1)$ and solving to obtain the parameters $\phi_1$ and $\theta_1$.
In the first case with $f=0.2$ the power is almost identical and in the second
case with $f=0.4$ the power is slightly higher for the ARMA$(1,1)$ approximation.
This suggests that long term memory in the fractional noise model has little
effect on the power when the length of the series is moderate,
as in this example with $n=50$ and $T=25$.
For sufficiently long time series, the effect on long memory is much more important
and the ARMA$(1,1)$ approximation does not hold.
\centerline{[Table 3 about here]}

\medskip
{\subsubti\noindent 2.5 Simulation Experiment}
\smallskip

The power function derived in eqn. (\ref{PowerFunctionTwoSided}) relies on the asymptotic normality of the maximum likelihood estimator and so it is helpful to
check its accuracy by simulation.
We do this by comparing the power function with the empirical power function, $\hat \Pi$.
For each simulated time series all parameters in the model were estimated by
exact maximum likelihood estimation and the $Z$-test was computed.
The empirical power, $\hat \Pi$, of a two-sided 5\% test
is then the proportion of times that the absolute value
of this $Z$-statistic exceeded 1.96 in absolute value and
the 95\% confidence interval for $\Pi$ is
$\hat \Pi \pm 1.96 \surd ( \hat \Pi (1-\hat \Pi )/N )$,
where $N$ is the number of simulations.
For each model and each parameter setting, $N=1,000$.

The model in eqn. (\ref{SimpleStepInterventionModel}) was simulated with
$n=50$ and $T=25$ and
$\AR(1)$ errors with $\phi_1=0,0.25,0.5,0.75$,
$\omega=\delta \sigma$, where
$\delta=0,\pm 0.25,..., \pm 2.0$.
The empirical power confidence limits and theoretical power
given by eqn. (\ref{PowerFunctionTwoSided}) are compared in
Figure \theFigureSimulatedPower.
It is seen that eqn. (\ref{PowerFunctionTwoSided})
provides an accurate approximation.
The $\IMA(1)$,
is a commonly occurring nonstationary time series model.
Figure \theFigureSimulatedPower\
compares the theoretical and empirical power for the case with $n=50$ and $T=25$
using a two-sided $Z$-test at the 5\% level.
Once again it is seen that eqn. (\ref{PowerFunctionTwoSided})
holds very well despite the small sample size.
The values selected for $\theta_1$ are positive since this is the most
common situation in practice.
The power improves, as expected, as $\theta_1$ increases from $0$ to $1$.
Notice that this model does not satisfy eqn. ({\ref{PierceCondition}).
The last column of Figure \theFigureSimulatedPower\
compares the empirical and theoretical power
in the case of fractionally differenced white noise, $p=q=d=0$ for
$f=0.0, 0.2, 0.3, 0.4$.
The approximation to the theoretical power improves with increasing $f$.
The simulations shown in Figure \theFigureSimulatedPower\ were repeated
using the likelihood-ratio test and essentially equivalent results were obtained.

\medskip
\centerline{[Figure \theFigureSimulatedPower\ about here]}
\medskip

In conclusion, the simulations in Figure \theFigureSimulatedPower\ suggest
that for practical purposes if $n$, $T$ and $n-T$ are not too small the
asymptotic theoretical power curve provides a good small sample approximation.
Alternatively, the simulations show that $\hat \omega$ is well approximated using
its large-sample approximation even for moderately small samples.
As already noted, $\sigma_{\hat \omega}$, must also be estimated by $\hat \sigma_{\hat \omega}$
using either the pre-intervention data or an estimate of its likely autocorrelation function.
In practice, as in the example in \S 2.1, a range of likely parameter values
are often used to indicate a range of possible power curves.

\medskip
{\subsubti\noindent 2.6 Model Uncertainty}
\smallskip

Box, Jenkins, and Reinsel (1994) found that both the ARMA$(1,1)$ and
IMA$(1)$ fit Series A, Chemical Process Concentrations about equally well.
Both models give similar one step ahead forecasts but the long run forecasts
are very different.
The situation is similar with the power functions for these two models.

Consider a hypothetical
step intervention which occurs immediately after the last observation.
In this case $T=198$ and the power curve as a function of $\omega$ is tabulated
for a few selected values in Table 7 for a two-sided 5\% test assuming
that $m$ post-intervention observations are available for $m=5$ and $m=50$.
When $m=5$ the power curves are quite similar but for $m=50$ the power
increases for the ARMA model but stays essentially the same in the case of
the IMA model.
For example, Table 7 shows that there is a 75\% chance of detecting a change of $0.6$ with just 5 post-intervention observations.

\centerline{[Table 4 about here]}

\medskip
{\subsubti\noindent 2.7 Forecast-Actuality Significance Test}
\smallskip

Box and Tiao (1976b) described an omnibus significance test for
detecting if an intervention has occurred.
If $a_t$, $t=T,...,n$ denote the one-step ahead prediction errors of
an assumed model, then the test statistic may be written,
$Q=\sum_{t=T}^n a_t^2/{\sigma_a^2}$.
If the intervention has no effect, $Q$ is approximately $\chi^2$-distributed on
$m=n-T+1$ df.
This significance test is easy to apply and does not require specification of an intervention
model and its estimation.
However, as might be expected, the loss of power can be considerable as will now be demonstrated.

As an example, consider the SIA model with a step intervention.
Then it can shown using eqn. (4) of Box and Tiao (1976b) that
$Q=||\omega 1_m^{\prime} \pi/\sigma_a +a/\sigma_a||^2$,
where $1_m$ denotes the $m$-dimension vector with $1$ in each position,
$a=(a_{T},...,a_{n})$, $\pi=(\pi_{i-j})$ is the lower triangular matrix
with $(i,j)$ entry $\pi_{i-j}$, where
$\pi_k$ is the coefficient of $B^k$ in the expansion
$\nabla^d \phi(B)/\theta(B)=1+\pi_1 B+\pi_2 B^2+...$.
So $Q$ has a $\chi^2$ distribution with $m$ df and noncentrality parameter
$\nu=(\omega^2/\sigma_a^2)||1_m^{\prime} \pi||^2$ and hence the large-sample
power function can be computed.
Figure \theFigureQPower\ compares the power of this significance test
with the SIA model hypothesis test for an example with $n=120$, $T=101$
and AR(1) errors.
Figure \theFigureQPower\
shows that the power of the significance test can be substantially less
than the intervention analysis hypothesis test.

\medskip
\centerline{[Figure \theFigureQPower\  about here]}

\bigskip
\centerline{\subti 3. ILLUSTRATIVE APPLICATIONS\/}

\medskip
{\subsubti\noindent 3.1 Traffic Safety and Public Policy}
\smallskip

On May 1, 1996, liquor bar closing time in Ontario was changed from 1 AM to 2 AM.
In a proposed intervention analysis we wished to examine the
possible effect of this change on late-night automobile fatalities.
The data for this study comprised the total number of fatalities
every month in Ontario during the hours of 11PM to 4AM for a period
of years before and after May 1, 1996.
For comparison we also collected similar time series data for
Michigan and New York State.
Data for this analysis were expensive to obtain since raw records
needed to be assembled, cleaned and aggregated from sources in various
jurisdictions.
Initially we planned to obtain monthly time series on the the
total number of fatalities from
January 1994 to December 1998.
This would yield $n=60$ observations and with the intervention
occurring at $T=36$.
At additional cost, we could obtain complete monthly time series
covering the period
January 1992 to December 1998
which corresponds to $n=84$ and $T=48$.
We were interested to know if $(n=60,T=36)$ or $(n=84,T=48)$ would be sufficient to
detect change of $\sigma$ or greater with a reasonably high probability,
where $\sigma$ is the standard deviation of the pre-intervention series.

Based on previous experience with similar time series
(Vingilis, et al., 1988)
we expected the time series will exhibit small autocorrelations which may be modelled by an $\AR(1)$ with parameter $\phi_1 \le 0.5$.
The intervention was expected to cause an increase in late-night fatalities, so a
one-sided upper-tail test is appropriate.
The power function in this case is
$\Pi(\delta)=$\ $1-\Phi(1.645-2.362\times \delta)$.
Table 5 shows the power of a 5\% upper-tail test for these two plans
for various $\phi_1$.
When $\phi_1=0.5$, Table 5 shows that $(n=84,T=48)$ has a 86.7\% chance of detecting a
step intervention whose magnitude is only one standard deviation of the error
component whereas the corresponding power for $(n=60,T=36)$ is 76.3\%.
The results of Table 5 demonstrated to our satisfaction and that of the granting agency,
that $(n=84,T=48)$ had a good chance of detecting a meaningful change
and was worth the extra expenditure.

\centerline{[Table 5 about here]}

\medskip
{\subsubti\noindent 3.2 Detecting Ozone Turnaround}
\smallskip

Tiao et al. (1990) used the SIA model with a ramp intervention with
$\AR(1)$ errors to model the trend in monthly deseasonalized stratospheric ozone and other environmental variables.
For simplicity Tiao et al. (1990) assumed that the mean of the pre-intervention series was known.
It may be shown that the expression obtained by Tiao et al. (1990, Appendix A) for
$\sigma_{\hat \omega}$
is exactly equal to
$\sigma_{\hat \omega} = 1/\surd {\cal I}_{2,2}$
using Table 1 with $n=T$ and $T=1$.
Table 6 compares this result
with the corresponding result obtained using the
exact expected Fisher information matrix given in eqn. (\ref{InformationMatrixGeneral})
for the same parameters as used in Tiao et al. (1990, Table 1).
When $\phi=0.8$, the difference is as high as 17\% but it decreases as the sample
size increases.
The approximation is very good for parameter values $0.6$ and less.
For most of the geophysical time series considered by Tiao et al. (1990) the degree
of autocorrelation is quite low, so this approximation works well.

\centerline{[Table 6 about here]}

Tiao et al. (1990, Table 2) also consider the number of years of monthly data
needed to detect a ramp intervention for several geophysical time series of interest.
In their computations it was assumed that $T=1$ and that
the mean was known.
Table 7 below computes the number of years of data needed for these time series
under the assumptions that the mean is unknown but that there are 30 years
of prior data.
The other assumptions about the data and the form of the intervention are
the same as in Tiao et al. (1990).
The parameter $\delta$ shown in the table was based on the information supplied
by Tiao et al. (1990).
Specifically,
$\delta=\omega/\left(12\times \hat \sigma \right)$
where $\hat \phi_1$ and $\hat \sigma$ are obtained from Tiao et al. (1990, Table 2)
and $\omega$ is obtained from Tiao et al. (1990, p.20,510).
Note that $\omega$ was divided by $12$ because the form
of the intervention used in Tiao et al. (1990) was $R_t^{(T+1)}/12$.
In conclusion, the estimate of the sample size required shown in Table 7
is in reasonable agreement with the results in Tiao et al. (1990).

\centerline{[Table 7 about here]}

\bigskip
\centerline{\subti 4. CONCLUDING REMARKS}
\medskip
We have shown how the power function for an intervention analysis may be computed
provided that we have an estimate of the $\ARIMA$ parameters in the
pre-intervention time series or in some closely related time series.
In the case of the SIA model with $\AR(1)$ or $\IMA(1)$ errors,
the power function can easily be computed using a hand calculator.
Such programs are freely available for the Texas Instruments TI-83 from the first author's
webpage.
{\it Mathematica\/} and S software for computing the power functions
and all tables and figures described in this paper
are also available there as well as various other supplements to this article.

The emphasis of this article has been on the use of the power function as an aid in selecting
the sample size.
In the case of the SIA model, if $\Pi(\omega^\prime)=1-\beta^\prime$ for a 5\% two-sided test
of ${\cal H}_0: \omega=0$ then
the usual 95\% confidence interval for $\omega$ will contain $0$ with probability $\beta^\prime$
when $\omega=\omega^\prime$.
So the power function may be used as an aid in choosing the sample size so that
a useful confidence interval is obtained.
Instead of the power function we could have focussed on the width of a suitable
interval estimate of $\omega$.
Since this also depends on an estimate of $\sigma_{\hat \omega}$ the methods presented
are applicable.
It may be noted that {\it overemphasis\/} on hypothesis tests has long been condemned
as was already noted many years ago by Cox (1977).
Nevertheless, as indicated by Cox (1977), such tests remain important in practice.

The power function depends strongly on the degree of autocorrelation in the
pre-intervention time series.
In the stratospheric ozone example, \S 2.2, a long pre-intervention series was
available which enabled the model to be accurately estimated.
In other cases, such as the traffic safety example, \S 2.1, the pre-intervention series
is either unavailable or quite short.
In such cases there may be prior information available which indicates a range
of likely models.
As discussed in \S 2.1, this may still be very useful for planning purposes.
A final note of caution, power computations should only be used before the analysis
of the data is done (Hoenig and Heisey, 2001; Lenth, 2001) and should never be used to compute the observed power after a test of hypothesis has already been carried out.

\bigskip
\centerline{\subti ACKNOWLEDGMENTS}
\parindent 0pt
This research was supported by grants from NIAAA and NSERC.
The authors would like to thank the Editor, an Associate Editor, two referees
and Dr. R.J. Kulperger for helpful suggestions.

\newpage
\centerline{\subti REFERENCES}
\parindent 0pt

\hind
Beran, J. (1994),
{\it Statistics for Long Memory Processes\/}.
London: Chapman and Hall.

\hind
Box, G.E.P., Hunter, W.G. and Hunter, J.S. (1978),
{\it Statistics for Experimenters\/},
New York: Wiley.

\hind
Box, G.E.P., Jenkins, G.M. and Reinsel, G.C.  (1994),
{\it Time Series Analysis: Forecasting and Control\/},
3rd Ed., San Francisco: Holden-Day.

\hind
Box, G.E.P. and Luce\~no, A. (1997),
{\it Statistical Control by Monitoring and Feedback Adjustment\/},
New York: Wiley.

\hind
Box, G.E.P. and Tiao, G.C. (1976a),
``Intervention Analysis with
Applications to Economic and Environmental Problems,''
{\it Journal of the American Statistical Association\/}, {\bf 70}, 70--79.

\hind
Box, G. E. P. and Tiao, G. C. (1976b),
``Comparison of Forecast and Actuality,''
{\it Applied Statistics\/} {\bf 25} (1976),  195--200.


\hind
Cao, J., Cleveland, W.S., Lin, D. and Sun, D.X. (2001),
``On the Nonstationarity of Internet Traffic,''
{\it Performance Evaluation Review: Proc. ACM Sigmetrics\/} {\bf 29}, 102-112.

\hind
Chang, I., Tiao, G.C. and Chen, C. (1988),
``Estimation of Time Series Parameters in the Presence of Outliers'',
{\it Technometrics\/} {\bf 30}, 193--204.



\hind
Cox, D.R. (1977),
``The Role of Statistical Signficance Tests,''
{\it Scand. J. Statist.\/} {\bf 4\/}, 49--70.

\hind
Golub, G. and Van Loan (1983),
{\it Matrix Computations\/},
Baltimore: John Hoptkins University Press.

\hind
Haslett, J. and Raftery, A. E. (1989),
``Space-time Modelling with Long-memory Dependence: Assessing Ireland's Wind Power Resource,''
{\it Applied Statistics\/} {\bf 38}, 1--21.

\hind Hipel, K.W. and McLeod, A.I. (1994). {\it Time Series
Modelling of Water Resources and Environmental Systems,\/}
Amesterdam: Elsevier.

\hind
Hoenig, J. M. and Heisey, D. M. (2001),
``The Abuse of Power:
The Pervasive Fallacy of Power Calculations for Data Analysis,''
{\it The American Statistician\/}, {\bf 55}, 19--24.

\hind
Jiang, W., Tsui, K.L. and Woodall, W.H. (2000),
``A New SPC Monitoring Method: The ARMA Chart,''
{\it Technometrics\/} {\bf 42}, 399--410.

\hind
Luce\~no, A. (1995),
``Choosing the EWMA Parameter in Engineering Process Control,''
{\it Journal of Quality Technology\/} {\bf 27}, 162--168.

\hind
K\"arner, O. (2002),
``On Nonstationarity and Antipersistency in Global Temperature Series,''
{\it Journal of Geophysical Research\/} {\bf 107} D20, 4415.

\hind
Lenth, R.V. (2001),
``Some Practical Guidelines for Effective Sample Size Determination,''
{\it The American Statistician\/} {\bf 55} 187--193.

\hind
Noakes, D. J. and Campbell, A. (1992),
''Use of Geoduck Clams to Indicate Changes in the Marine Environments of Ladysmith Harbour, British Columbia,''
{\it EnvironMetrics\/} {\bf 3}, 81--97.

\hind
Pierce, D.A. (1972),
``Least Squares Estimation in Dynamic-disturbance Time Series Models,''
{\it Biometrika}  {\bf59\/}, 73--78.

\hind
Reinsel, G. C. (2002),
``Trend Analysis of Upper Stratospheric Umkehr Ozone Data for Evidence of Turnaround,''
{\it Geophysical Research Letters\/} {\bf 29(10)}, doi:10.1029/2002GL014716.

\hind
Reinsel, G. C., Weatherhead, E. C., Tiao, G. C., Miller, A. J., Nagatani, R. M., Wuebbles, D. J., and Flynn, L. E. (2002),
``On Detection of Turnaround and Recovery in Trend for Ozone,''
{\it Journal of Geophysical Research\/} {\bf 107} (D10), doi:10.1029/2001JD000500.

\hind
Roberts, (1991),
{\it Data Analysis for Managers with Minitab.\/}
2nd Ed. San Francisco: The Scientific Press.

\hind
Tiao, G.C., Reinsel, G.C., Xu, D., Pedrick, J.H., Zhu, X., Miller, A.J., DeLuisi, J.J.,
Mateer, C.L. and Wuebbles, D.J. (1990),
``Effects of Autocorrelation and Temporal Sampling Schemes on Estimation
of Trend and Spatial Correlation,''
{\it Journal of Geophysical Research\/} {\bf 95} D12, 20,507--20,517.

\hind
Weatherhead, E.C., Reinsel, G.C., Tiao, G.C., Meng, X.L., Choi, D., Cheang, W.K.,
Keller, T., DeLuisi, J., Wuebbles, D.J., Kerr, J.B., Miller, A.J., Oltmans, S.J.
and Frederick, J.E. (1998),
``Factors Affecting the Detection of Trends: Statistical Considerations
and Applications to Environmental Data'',
{\it Journal of Geophysical Research\/} {\bf 103} D14, 17,149--17,161.

\hind
Vingilis, E., Blefgen, H., Lei, H., Sykora, K. and Mann, R. (1988).
``An Evaluation of the Deterrent Impact of Ontario's 12-hour Licence Suspension Law,''
{\it Accident Analysis and Prevention\/} {\bf 20}, 9--17.

\hind
Wolfram, S. (1999),
{\it The Mathematica Book\/},
4th Ed., Wolfram Media/Cambridge University Press, Champaign/Cambridge.



\vfill\eject
\clearpage
\hoffset=0truecm
\begin{figure}
\scalebox{1.0}{\includegraphics{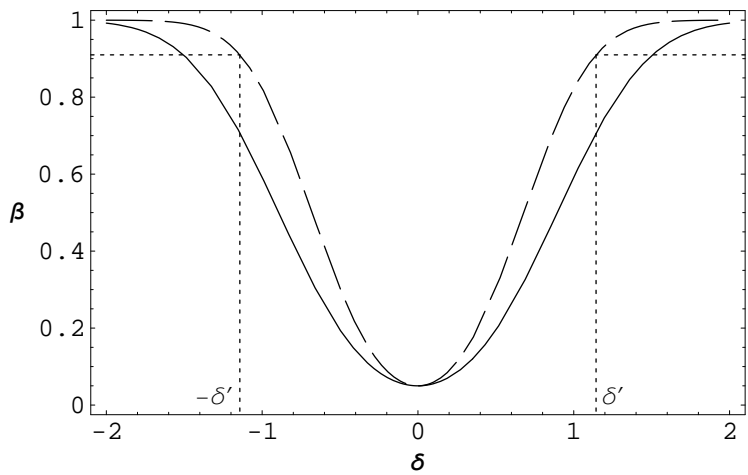}}
\caption{
\tti Comparison of Power Curves For $n=50, T=25$ and $n=\infty, T=25$.
The solid curve shows for $n=50,T=25$
and
the dashed curve, $n=\infty, T=25$.
The approximate detection limit, $\delta^{\prime}\doteq 1.143$ is also
shown.
}
\end{figure}

\vfill\eject

\clearpage
\hoffset=0truecm
\voffset=0truecm
\begin{figure}
\scalebox{0.95}{\includegraphics[3cm,7cm][5cm,7cm]{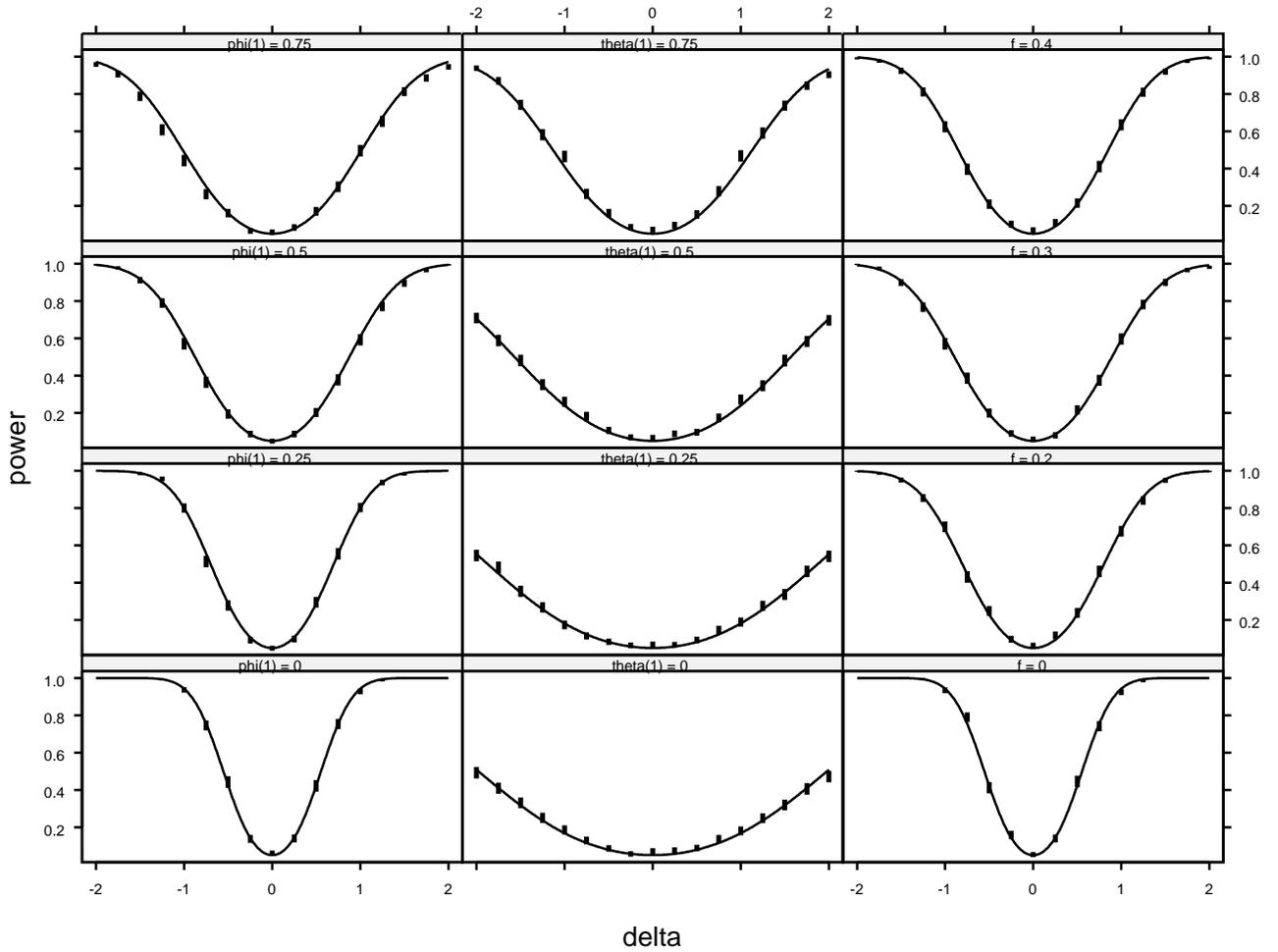}}
\vskip2in
\caption{
\tti Comparison of Empirical and Theoretical Asymptotic Power
in the SIA Model with AR(1), IMA(1) and Fractionally-Differenced White Noise.
The parameter $\delta=\omega/\sigma$ is the rescaled step size.
The solid curve shows the theoretical power defined in eqn. (\ref{PowerFunctionTwoSided}).
The vertical bars show the width of a 95\% confidence
interval for the empirical power in 1,000 simulations of the model.
The AR(1) and IMA(1) parameters $\phi_1$ and $\theta_1$ are denoted by
phi(1) and theta(1) in the diagram.
}
\end{figure}

\vfill\eject

\clearpage
\hoffset=0truecm
\begin{figure}
\scalebox{1.}{\includegraphics[3cm,3cm][2cm,2cm]{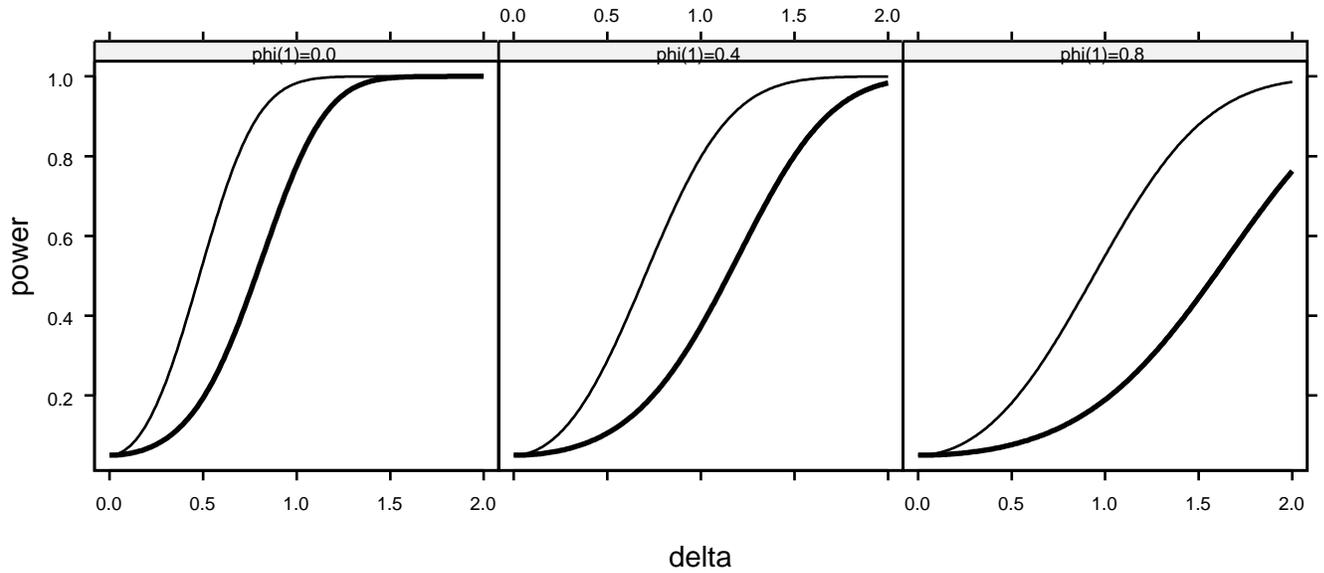}}
\caption{
\tti Comparison of Power Functions for a SIA Model with a Step Intervention
with AR(1) Errors
and the Forecast-Actuality Significance Test For a Two-Sided Test
at the 5\% Level.
The model parameters are $n=120$, $T=101$,
${\tti delta\/}\ = \delta=\omega/\sigma$
and ${\tti phi(1)\/}\ = \phi_1$.
The solid thin curve shows the SIA Model based hypothesis test and
the solid thick curve shows the omnibus significance test using $Q$.
Since both power functions are symmetric about $\delta=0$ only the
upper half is shown.
}
\end{figure}

\vfill\eject

\clearpage
\hoffset=0truecm
\begin{table}
\caption{\tti Information Matrix for Simple Intervention Analysis
with AR(1) Errors.
The table gives the $(1,2)$ and $(2,2)$ entries,
${\cal I}_{1,2}/\sigma_a^2$ and ${\cal I}_{2,2}/\sigma_a^2$.
For each intervention type, ${\cal I}_{1,1}/\sigma_a^2=n (1-\phi_1)^2$
and the $(2,1)$ entry is obtained by symmetry.
}
$$\vbox{\halign{\strut
 \hfill#\quad    
&\hfill#\quad    
&#\hfill   
&#\hfill   
\cr             
\noalign{\hrule}
\noalign{\smallskip}
\noalign{\hrule}
\noalign{\smallskip}
Type& & Information Matrix Entries
\cr
\noalign{\smallskip}
\noalign{\hrule}
\noalign{\smallskip}
Step & ${\cal I}_{1,2}/\sigma_a^2=$ & $ (n-T)(1-\phi_1)^2+1-\phi_1$ \cr
& ${\cal I}_{2,2}/\sigma_a^2=$ & $ (n-T)(1-\phi_1)^2+1$
\cr
\noalign{\smallskip}
\noalign{\hrule}
\noalign{\smallskip}
Pulse & ${\cal I}_{1,2}/\sigma_a^2=$ & $ 1-\phi_1^2 $ \cr
& ${\cal I}_{2,2}/\sigma_a^2=$ & $ 1-\phi_1^2 $
\cr
\noalign{\smallskip}
\noalign{\hrule}
\noalign{\smallskip}
Ramp & ${\cal I}_{1,2}/\sigma_a^2=$ & $ \left( 1 + n - T \right) \,\left( 1 - \phi_1  \right) \,
    \left( 2 + n - T - \left( n - T \right) \phi_1  \right) /2 $  \cr
& ${\cal I}_{2,2}/\sigma_a^2=$ & $\left( 1 + n - T \right) \,
    ( 6 + 7\,n + 2\,n^2 - 7\,T - 4\,n\,T + 2\,T^2 - 8\,n\,\phi_1 $ \cr
& &$-4\,n^2\,\phi_1  + 8\,T\,\phi_1  + 8\,n\,T\,\phi_1  - 4\,T^2\,\phi_1  +
      n\,{\phi_1 }^2 + 2\,n^2\,{\phi_1 }^2 - T\,{\phi_1 }^2 $ \cr
& &$-4\,n\,T\,{\phi_1 }^2 + 2\,T^2\,{\phi_1 }^2  )/{6}$
\cr
\noalign{\smallskip}
\noalign{\hrule} \cr
}}$$
\end{table}

\begin{table}
\caption{\tti Information Matrix for Simple Intervention Analysis
with $\IMA(1)$ Errors. For $\theta_1=0$ set $\theta_1^0=1$.
The table gives the $(1,2)$ and $(2,2)$ entries,
${\cal I}_{1,2}/\sigma_a^2$ and ${\cal I}_{2,2}/\sigma_a^2$.
For each intervention type, ${\cal I}_{1,1}/\sigma_a^2=(n-1) / (1-\theta_1)^2$
and the $(2,1)$ entry is obtained by symmetry.
}
$$\vbox{\halign{\strut
 \hfill#\quad    
&\hfill#\quad    
&#\hfill   
&#\hfill   
\cr             
\noalign{\hrule}
\noalign{\smallskip}
\noalign{\hrule}
\noalign{\smallskip}
Type& & Information Matrix Entries
\cr
\noalign{\smallskip}
\noalign{\hrule}
\noalign{\smallskip}
Step & ${\cal I}_{1,2}/\sigma_a^2=$ & $ (1-\theta_1)^{-2} (1-\theta^{n+1-T}) $ \cr
& ${\cal I}_{2,2}/\sigma_a^2=$ & $ (1-\theta_1^2)^{-1} (1-\theta^{2(n+1-T)}) $
\cr
\noalign{\smallskip}
\noalign{\hrule}
\noalign{\smallskip}
Pulse & ${\cal I}_{1,2}/\sigma_a^2=$& $ (1-\theta_1)^{-1} \theta_1^{n-T} $ \cr
& ${\cal I}_{2,2}/\sigma_a^2=$& $(1+\theta_1)^{-1} 2  (1+\theta^{2 (n -T) + 1}) $
\cr
\noalign{\smallskip}
\noalign{\hrule}
\noalign{\smallskip}
Ramp &${\cal I}_{1,2}/\sigma_a^2=$ & $ (1-\theta_1)^{-3} (n+1-T+\theta_1^{n+2-T}-(n+2-T)\theta) $ \cr
&${\cal I}_{2,2}/\sigma_a^2=$ & $ (1+\theta)^{-1} (1-\theta_1)^{-3} (  2\,{\theta }^{2 + n + T}\, ( 1 + \theta  )  -{\theta }^{4 + 2\,n}$ \cr
 &&$ + {\theta }^{2\,T}\,( n+1 - T - 2\,\theta  -( 2 + n - T ) \,{\theta }^2 )  )  $
\cr
\noalign{\smallskip}
\noalign{\hrule} \cr
}}$$
\end{table}

\vfill\eject
\clearpage

\begin{table}
\caption{\tti Power Function, $\Pi(\delta)$, for Fractionally Differenced White Noise
With Parameter $f$ and The Approximating ARMA$(1,1)$ Model
for a Two-sided 5\% Level Test in SIA Step Intervention Model
with $n=50$ and $T=25$.
The first entry in each pair is for
the fractional model and the second the ARMA$(1,1)$ model.
The parameters in the approximating ARMA model are respectively
$\phi_1=0.667, \phi_2=0.451$
and
$\phi_1=0.875, \phi_2=0.405$
corresponding respectively to $f=0.2$ and $f=0.4$.
}
$$\vbox{\halign{\strut
 #\hfill \quad    
&\hfill#\quad    
&\hfill#\quad    
&\hfill#\quad    
&\hfill#\quad    
&\hfill#\quad    
&\hfill#\quad    
\cr             
\noalign{\hrule}
\noalign{\smallskip}
\noalign{\hrule}
\noalign{\smallskip}
$\delta$\qquad&$f=0.2$\qquad&$f=0.4$\qquad
\cr
\noalign{\smallskip}
\noalign{\hrule}
\noalign{\smallskip}
$0$&$0.050,0.050$&$0.050,0.050$
\cr
$0.5$&$0.198,0.202$&$0.086,0.076$
\cr
$1.$&$0.602,0.612$&$0.198,0.156$
\cr
$1.5$&$0.914,0.920$&$0.384,0.291$
\cr
$2.$&$0.993,0.994$&$0.602,0.468$
\cr
$2.5$&$1.000,1.000$&$0.792,0.651$
\cr
$3.$&$1.000,1.000$&$0.914,0.805$
\cr
\noalign{\smallskip}
\noalign{\hrule}
}}$$
\end{table}

\vfill\eject
\clearpage
\begin{table}
\caption{\tti Power Comparison for Step Interventions with ARMA(1,1) and IMA(1) Errors for Series A with $n=197+m$
and $T=198$.  The models' other parameters are respectively,
\{$\phi_1=0.9087, \theta_1=0.5758, \sigma_a=0.3125$\}
and
\{$\theta_1=0.7031,  \sigma_a=0.3172$\}.
}
\bigskip

\begin{center}
\begin{tabular}{ccccc}
&\multicolumn{2}{c} \hfill ARMA$(1,1)$\hfill &\multicolumn{2}{c}\hfill IMA$(1)$\hfill\cr
\noalign{\smallskip}
\noalign{\hrule}
$\omega$&$m=5$&$m=50$&$m=5$&$m=50$
\cr
\noalign{\smallskip}
\noalign{\hrule}
\noalign{\smallskip}$0.2$&$0.141$&$0.205$&$0.141$&$0.143$
\cr
$0.3$&$0.258$&$0.398$&$0.258$&$0.264$
\cr
$0.4$&$0.415$&$0.621$&$0.416$&$0.425$
\cr
$0.5$&$0.588$&$0.809$&$0.589$&$0.600$
\cr
$0.6$&$0.745$&$0.925$&$0.746$&$0.756$
\cr
$0.7$&$0.863$&$0.978$&$0.864$&$0.872$
\cr
\end{tabular}
\end{center}

\end{table}

\clearpage

\begin{table}
\caption{\tti Power Comparison for AR(1) Errors for $(n=60,T=36)$ and $(n=84,T=48)$.
The first entry in each column corresponds to $(n=60,T=36)$ and the second $(n=84,T=48)$.
}

$$\vbox{\halign{\strut
 \hfill#\quad    
&\hfill#\quad    
&\hfill#\quad    
&\hfill#\quad    
&\hfill#\quad    
&\hfill#\quad    
&\hfill#\quad    
&\hfill#\quad    
&\hfill#\quad    
\cr             
\noalign{\hrule}
\noalign{\smallskip}
\noalign{\hrule}
\noalign{\smallskip}
$\delta$&$\phi_1=0$&$\phi_1=0.25$&$\phi_1=0.5$&$\phi_1=0.75$
\cr
\noalign{\smallskip}
\noalign{\hrule}
\noalign{\smallskip}
$0.000$&$0.050,0.050$&$0.050,0.050$&$0.050,0.050$&$0.050,0.050$
\cr
$0.250$&$0.245,0.306$&$0.186,0.226$&$0.146,0.170$&$0.124,0.135$
\cr
$0.500$&$0.604,0.736$&$0.444,0.555$&$0.321,0.395$&$0.253,0.288$
\cr
$0.750$&$0.889,0.961$&$0.729,0.848$&$0.550,0.664$&$0.431,0.493$
\cr
$1.000$&$0.985,0.998$&$0.914,0.973$&$0.763,0.867$&$0.624,0.700$
\cr
$1.250$&$0.999,1.000$&$0.983,0.998$&$0.904,0.964$&$0.790,0.857$
\cr
$1.500$&$1.000,1.000$&$0.998,1.000$&$0.971,0.994$&$0.903,0.946$
\cr
$1.750$&$1.000,1.000$&$1.000,1.000$&$0.994,0.999$&$0.963,0.984$
\cr
$2.000$&$1.000,1.000$&$1.000,1.000$&$0.999,1.000$&$0.989,0.996$
\cr
\noalign{\smallskip}
\noalign{\hrule}
}}$$
\end{table}

\vfill\eject
\clearpage

\begin{table}
\caption{\tti Comparison of Exact and Approximate Methods.
The function $g(T,\phi)$ defined in Tiao et al. (1990) was computed
using exact form of the information matrix eqn. (\ref{InformationMatrixGeneral})
and the approximation eqn. (\ref{InformationMatrixPierce}) for selected
parameter values given in Table 1 of Tiao et al. (1990).
The entries in the table show the percentage difference,
$100\times ({\rm EXACT} - {\rm APPROXIMATE})/{\rm EXACT}$.
}

$$\vbox{\halign{\strut
 \hfill#\quad    
&\hfill#\quad    
&\hfill#\quad    
&\hfill#\quad    
&\hfill#\quad    
\cr             
\noalign{\hrule}
\noalign{\smallskip}
\noalign{\hrule}
\noalign{\smallskip}
Number&& \cr
of&$\phi=0.6$&$\phi=0.8$ \cr
Years&& \cr
\noalign{\smallskip}
\noalign{\hrule}
\noalign{\smallskip}$6$&$-6$&$-17$
\cr
$7$&$-5$&$-15$
\cr
$8$&$-5$&$-13$
\cr
$9$&$-4$&$-11$
\cr
$10$&$-4$&$-10$
\cr
\noalign{\smallskip}
\noalign{\hrule}
}}$$
\end{table}

\vfill\eject
\clearpage

\begin{table}
\caption{\tti Number of Years, $n^{*}$, For 90\% Probability of Detecting a Prescribed
Trend, $\delta$ Using a Two-Sided $5$\% Test Given 30 Years of Prior Data
And Assuming $\AR(1)$ Errors With Estimated Parameter $\hat \phi_1$.
The last line of the table shows the comparable values given in
Tiao et al. (1990, Table 2).
}

$$\vbox{\halign{\strut
#\hfill\qquad    
&#\hfill\qquad   
&#\hfill\qquad   
&#\hfill\qquad   
&#\hfill\qquad   
&#\hfill\qquad   
\cr             
\noalign{\hrule}
\noalign{\smallskip}
\noalign{\hrule}
\noalign{\smallskip}
&Tateno&Hohen.&Wakkan&Bulawayo&Abidajan\cr
\noalign{\smallskip}
\noalign{\hrule}
\noalign{\smallskip}
$\hat \phi_1$&$0.32$&$0.05$&$0.14$&$0.43$&$0.65$
\cr
\noalign{\smallskip}
$\omega$&$0.003$&$0.003$&$0.2$&$0.2$&$0.2$
\cr
\noalign{\smallskip}
$\delta$&$0.00758$&$0.00543$&$0.01042$&$0.01282$&$0.01111$
\cr
\noalign{\smallskip}
$n^{*}$&$11.6$&$12.1$&$8.0$&$8.6$&$12.0$
\cr
\noalign{\smallskip}
$n^{*}_{\rm Tiao}$&$14$&$14$&$10$&$10$&$13$
\cr
\noalign{\smallskip}
\noalign{\smallskip}
\noalign{\hrule}
}}$$
\end{table}

\vfill\eject
\clearpage

\begin{table}
\caption{\tti Power Comparisons of Dynamic Step Intervention Model with
Simple Step Intervention when $n=50$ and $T=25$.
The first entry in each triplet shows the theoretical power of a 5\% two-sided test of ${\cal H}_0: g=0$
where $g=\omega_0^{(1)}/(1-\delta_1)$
in the dynamic step intervention model
$z_t = \xi+ \omega_0^{(1)}/(1-\delta_1 B)  S_t^{(T)} +  a_t/(1-\phi_1 B)$ with
$\xi=0, \phi_1 = 0.5$ and $\sigma_a^2=1$.
The second entry is the theoretical power of a 5\% test of
${\cal H}_0: \omega_0^{(2)}=0$
in the SIA model,
$z_t = \xi+ \omega_0^{(2)} S_t^{(T)} +  a_t/(1-\phi_1 B)$,
where $\omega_0^{(2)}=\omega_0^{(1)}/(1-\delta_1)$ and all other parameters are the same
as in the dynamic model.
The third entry is the empirical power, based on 1000 simulations, for a two-sided 5\% test of
${\cal H}_0: \omega_0^{(2)}=0$
when the SIA model is fitted to a time series generated by the
dynamic step intervention model.
}

$$\vbox{\halign{\strut
 \hfill#\quad    
&\hfill#\quad    
&\hfill#\quad    
&\hfill#\quad    
\cr             
\noalign{\hrule}
\noalign{\smallskip}
\noalign{\hrule}
\noalign{\smallskip}
$\delta_0$&$\omega_0=0.5$&$\omega_0=0.75$&$\omega_0=1.0$
\cr
\noalign{\smallskip}
\noalign{\hrule}
\noalign{\smallskip}$0.25$&$0.226,0.252,0.241$&$0.416,0.490,0.466$&$0.879,0.972,0.880$
\cr
$0.50$&$0.439,0.490,0.445$&$0.745,0.827,0.758$&$0.997,1.000,0.974$
\cr
$0.75$&$0.673,0.732,0.692$&$0.937,0.972,0.932$&$1.000,1.000,0.955$
\cr
\noalign{\smallskip}
\noalign{\hrule}
}}$$

\end{table}

\vfill\eject
\clearpage
\centerline{\subti Appendix A: Derivation of the Information Matrix}
\bigskip
The loglikelihood function, apart from a constant, may be written,
\begin{equation}
L(\lambda_1, \lambda_2, \sigma_a^2)=
-\log(\sigma)-\log(\det(\Gamma_n))-{1\over 2 \sigma_a^2} y^{\prime} \Gamma_n^{-1} y,
\label{loglikelihoodA}
\end{equation}
where $y$ is the column vector of length $n-d$ with $t$-th entry
$\nabla^d z_t -\xi -\omega \nabla^d S_t^{(T)},\ t=d+1,\ldots,n$.
Then $\partial y / \partial \xi= (-1,\ldots,-1)$.
Similarly $\partial y / \partial \omega=(-S_1^{(T)},\ldots,-S_n^{(T)})$.
Hence,
\begin{eqnarray}
{\cal I}_{\lambda_1} & =
&-E(\partial^2_{\lambda_1,\lambda_1} L(\lambda_1, \lambda_2, \sigma_a^2)) \cr
&=& {1\over  \sigma_a^2} J^{\prime} \Gamma_n^{-1} J,
\label{DerivationInformationMatrixA}
\end{eqnarray}
where $J$ is as in eqn. (\ref{InformationMatrixGeneral}).
Since
$E(\partial^2 L(\lambda_1, \lambda_2, \sigma_a^2) /(\partial\lambda_1 \partial\lambda_2))=0$
and
$E(\partial^2 L(\lambda_1, \lambda_2, \sigma_a^2) /(\partial\lambda_1 \partial\lambda_2))=0$,
the information matrix is block diagonal.

\vfill\eject
\clearpage
\centerline{\subti Appendix B: Interventions With A Dynamic Response}
\bigskip

For completeness we also discuss
the intervention analysis model with a dynamic response to the intervention
which may be written,
\begin{equation}
\nabla^d  z_t =
\xi+
\omega(B)/\delta(B) \nabla^d B^b I_t^{(T)} +
\nabla^{-f} \frac{\theta(B)}{\phi(B)} a_t,\quad t=1,\ldots,n,
\label{DynamicInterventionModel}
\end{equation}
where
$\omega(B)=\omega_0 + \omega_1 B + \ldots \omega_r B^r$
and
$\delta(B)=\delta_0 - \delta_1 B - \ldots \delta_s B^s$.
For stability of the transfer function it is assumed that all roots of $\delta(B)=0$
lie outside the unit circle.
As in Appendix A, the exact information matrix for the
parameters $\lambda_1 = (\xi,\omega_0,\ldots,\omega_r,\delta_1,\ldots,\delta_s)$
${\cal I}_{\lambda_1}=\sigma_a^{-2} J^{\prime} \Gamma_n^{-1} J$
where
$J$ is an $n-d \times (2+r+s)$ matrix
with rows
$(1$, $u_t, \ldots, u_{t-r}$, $v_t, \ldots, v_{t-s})$
for $t=1,\ldots,n-d$,
where
$u_{t-j} = \nabla^d \left(1/\delta(B)\right) I_{t-j}^{(T)}$
and
$v_{t-j} = \nabla^d \left(\omega(B)/\delta(B)\right) I_{t-j}^{(T)}$.
Alternatively the large-sample approximation given in Pierce (1972) may
be used.
The steady-state gain (Box, Jenkins and Reinsel, 1994, \S 10.1.1),
which measures the long-run change of the intervention,
is defined by
$g=$\ $(\omega_0 + \ldots + \omega_r)$$/(1-\delta_1 - \ldots -\delta_s)$.
The maximum likelihood estimates for the model may be used to form the estimate
of $g$, $\hat g$.
Using a Taylor series linearization, the standard deviation of $\hat g$ is
given by
$\sigma_{\hat g}= \surd (d_{\zeta}^{\prime} V_{\zeta} d_{\zeta})$,
where $V_{\zeta}$ is obtained by dropping the first row and column from
${\cal I}_{\lambda_1}^{-1}$ and
$d_{\zeta}=$\
$(\partial g / \partial \omega_0,\ldots,\partial g / \partial \omega_r$,
$\partial g / \partial \delta_1,\ldots,\partial g / \partial \delta_s)$.
For dynamic intervention analysis models we may consider testing
${\cal H}_0: g = 0$ using the $Z$ test.
Notice that, when $s>0$ we need estimates of all
parameters in the full intervention model to estimate $\sigma_{\hat g}$.
This limits the applicability of this approach since even if the pre-intervention
series is known, it is not likely that such precise information is available for the
intervention parameters.
Often the SIA model can be used to get an approximation to the power in this case.

As a numerical illustration, consider the dynamic step intervention model,
$z_t = \xi+ \omega_0{(1)}/(1-\delta_1 B)  S_t^{(T)} +  a_t/(1-\phi_1 B),
t=1,\ldots,n$.
Taking $n=50, T=25, \xi=0, \phi_1 = 0.5$ and $\sigma_a^2=1$,
Table 8 below compares the power of a 5\% two-sided test ${\cal H}_0: g = 0$,
where $g=\omega_0{(1)}$,
with that of the Z-test ${\cal H}_0: \omega_0^{(2)} = 0$
in the corresponding SIA model defined by
$z_t = \xi+ \omega_0^{(2)}  S_t^{(T)} +  a_t/(1-\phi_1 B)$
where $\omega_0^{(2)}=g$ and the other parameter settings are the same.
On an intuitive basis, the effect
in the SIA model is slightly
larger so one might expect the power in the SIA model to be
slightly larger.
Table 8 shows, comparing the first two entries in each triplet,
that this is exactly what happens.
The third entry in each triplet in Table 8 is the empirical power of
a two-sided 5\% test of ${\cal H}_0: \omega_0^{(2)}=0$
when the SIA model is fitted to a time series generated by the
dynamic step intervention model.
One thousand simulations were used for each model.
The empirical power is predicted well by the theoretical asymptotic
power for the SIA model.
These simulations were repeated with various values of the parameter
$\phi$ and similar results where found when $-1 < \phi \le 0.5$.
For $\phi_1 > 0.5$, there was a much bigger difference between
the asymptotic theoretical power of the dynamic and step models.
For example with $\phi_1=0.9$, $\omega_1=0.75$ and $\delta_1 = 0.75$, the asymptotic
power for the two-sided 5\% level gains test was only 0.199 whereas the
predicted power using a SIA step intervention was 0.972.
The empirical power of the two-sided 5\% level test of ${\cal H}_0: \omega_1 = 0$
in the step SIA model was 0.283.
The general conclusion reached was that the step SIA model provides a useful
approximation to the more complicated dynamic step intervention model provided
the autocorrelation is not too large.
Further simulation results are available in the online supplements.

\end{document}